# A useful ethics framework for mathematics teachers

Lucy Rycroft-Smith[1], Dennis Müller[2], Maurice Chiodo[3], Darren Macey[4]

**Abstract:** An argument is presented that mathematics teachers would benefit from an ethics framework that allows them to reflect upon their position regarding the teaching of ethics within mathematics. We present such a framework, adapted from Chiodo and Bursill-Hall's (2018) levels of ethical engagement for mathematicians, with discussion of its applications and implications. Awareness of ethical and social responsibility amongst mathematicians and mathematics educators is a newly established field, and yet it is crucially important, if we are to use mathematics for good and the avoidance of harm – in other words to use mathematics for justice, belonging, peace, care, education, lawmaking, building, creativity, equity, connection, fairness, rulemaking and order. Here, we present a way to extend ideas from ethics in the field of *mathematics for mathematicians* to the field of *mathematics education for mathematics teachers*, aiming to contribute to this transdisciplinary work and to advance the field further. We also explore links to ideas of equity, and silence and voice, drawing on the consequences of a deeper examination of the entanglements of mathematical-ethical dilemmas for mathematics teachers both inside and outside the classroom.

**Key words:** ethics in mathematics, mathematics education, mathematics for social justice, teaching framework

# Introduction

Mathematics itself has been unhelpfully conceived as pure, value-free, apolitical, and neutral – and therefore not ground where ethical work is even considered necessary (e.g., Ernest 2016a; 2020), which is one of the first barriers to a consideration of ethics in mathematics. Mathematics curricula, resources and assessment generally do not include useful information on ethical awareness, although they do often attend to matters of context. However, ethical issues exist in the mathematics classroom – of course they do – not only in the materials students interact with but also in the relationships, choices and agencies exercised in mathematical, and indeed any pedagogical, activities. Mathematics teachers deal in ethics every day, and the decision to remove these ideas from the sanitised version of 'school mathematics' we are given in curricula and associated resources is a form of abstraction that has, we argue, rendered invisible information that is highly salient and crucially important. This point is not new. Many people for many years have argued that school mathematics is artificial, having been contextwashed in some way. For example, in the context of the 1902 French educational reforms the mathematician Emile Borel argued for ways "to introduce more life and a greater sense of reality to our mathematics education" (Borel 1904, p. 436, translation by and quoted in Gispert & Schubring, 2011, p. 78). Where ideas with ethical dimensions do exist in formal materials, these tend to be constrained to areas broadly construed as statistics and data handling. This, along with other faux fault lines, can be seen as creating a false distinction between mathematics and statistics education where truth and fudging, certainty and uncertainty, and purity and interpretation are set up in false dichotomies, deluding students and teachers alike that dealing with these two types of mathematical beasts require different skills, methods and sensibilities.

---

[1] University of Cambridge, United Kingdom. Email: lr444@cam.ac.uk
[2] RWTH Aachen University, Aachen, Germany. Email: dennis.mueller3@rwth-aachen.de
[3] Department of Pure Mathematics and Mathematical Statistics, University of Cambridge. King's College, Cambridge. United Kingdom. Email: mcc56@cam.ac.uk
[4] University of Cambridge, United Kingdom. Email: djm249@cam.ac.uk

Here, we propose that the activity of mathematical modelling – best described as using mathematics to explain, define, estimate and test real problems and solutions, often in classes or systems (e.g., Arseven, 2015) – allows us to see beyond these naïve binaries and focus on the real, useful activities inherent in "reading the world with mathematics" (Gutstein, 2003, p. 45). We argue that mathematical modelling is not only central to all of school mathematics, but, however construed or constituted, has ethics embedded deep within its structures. Hence, we see the need for teachers of school mathematics at all levels to be supported by a framework for the development of their ethical awareness and responsibility so that they are able to engage with these issues in a way that, in turn, supports the development of ethical mathematicians in their classrooms and also supports equitable access to mathematics. This is emphatically not to say that only students who go on to be academic mathematicians, or even engineers, economists or statisticians need this guidance, but that all those who practice mathematics at any level deserve to understand how to use it as a tool for justice and belonging.

Chiodo and Bursill-Hall (2018) propose a useful level-based system for mathematicians' ethical engagement or consciousness, suggesting these levels are neither discrete nor linear exactly, but a useful observation from which to begin more nuanced discussions of this important topic that will broaden and deepen over time.

- Level -1: Actively obstructing efforts to address ethics in mathematics (a new level not explicitly stated in Chiodo and Bursill-Hall (2018), but suggested here as an extension)
- Level 0: Believing there is no ethics in mathematics
- Level 1: Realising there are ethical issues inherent in mathematics
- Level 2: Doing something: speaking out to other mathematicians
- Level 3: Taking a seat at the tables of power
- Level 4: Calling out the bad mathematics of others

In this chapter, we update and extend this theoretical observation to a framework for use in the context of mathematics education, in particular applying it to mathematics teachers of students of all ages. The core message remains unchanged – that the ethical mathematics teacher, like the ethical mathematician, seeks to not just do good, but, above all else, tries to prevent harm (Müller, 2022). However, here it is interpreted as, for example, avoiding the constituting of mathematics education as symbolic violence (Archer et al., 2018). Further, the richer sense of these proposed levels for the mathematics education context is that the ethical mathematics teacher has a parallel responsibility: to develop a sense of ethics in the mathematicians they teach, too. Here we imagine what that may look like in practice, drawing on the literature around ethics in mathematics and mathematics education more widely in order to provide such a framework for mathematics educators. We also consider the ways in which statistics education may provide a blueprint in this area that mathematics education can draw upon, with more embedded arguments for incorporating ethical responsibility in learning and teaching throughout students' education (Ridgway, 2016).

# Background

## Conceptions of ethics

Most definitions of ethics are concerned with decision-making in the pursuit of doing good and the avoidance of harm. Several researchers, such as Chiodo and Clifton (2019), Ernest (2018), and Müller et al. (2022), have suggested with rich examples that mathematics can be used to do harm as well as good, albeit often unintentionally. Similarly, Ernest (2016a) raises the question of whether mathematics education, at least in the White Western tradition, can claim to be in any such way 'ethical' when it currently negatively impacts such a large number of students. Over the years, Hippocratic oaths have been proposed for mathematicians (for an overview see Müller et al., 2022) and for teachers (Cody, 2007), drawing on ideas used in medicine to create ethical frameworks that consider ideas of autonomy, non-exploitation, privacy, consideration of intentional and unintentional consequence, need for collaboration, bias and fairness. But just as oaths alone may be insufficient to solve the problem (Müller et al., 2022), so would be a naïve translation of abstract philosophical concepts. Classrooms are messy, fast-paced, rich with interlaced decisions, and thick with power moves, and any ethics framework for teachers should, we believe, start here, and be first and foremost *practical*.

We draw on Danks (2021), who proposes the useful concept of *translational ethics*, applied in this context to digital technologies, from principles to practice and grounded in actual practices. They outline several types of oversimplification of ethics: as mere personal belief, as relativism with no right or wrong and further no better or worse, as compliance with law, as constraint only (negative only), and as freedom only (positive only). By contrast, they conceive of ethics as embedded in the practical, complex, multidimensional, almost always associated with several possible acceptable actions or outcomes, and most importantly concerned with *values* – asking "what values ought we to have?" and "given our values, how ought we act?"(ibid, p. 2).

## Values

Values are defined in mathematics education by Bishop et al. (2001, p. 1) as "those deep, affective qualities which education aims to foster through the school subject of mathematics and are a crucial component of the classroom affective environment." They are often embedded, implicit, and can be entirely invisible to the untrained eye. In abstract settings, such as mathematics, they can be especially hard to detect. Moreover, some people may deliberately choose to use abstraction as a way to avoid 'subjective' decision-making, and therefore ethical reasoning. Additionally, Chiodo and Vyas (2019, pp. 2-3) argue for the need of ethics in mathematics education because "abstract pure mathematics courses are often taught by individuals with little or no experience working or engaging with structures outside academia, thus preventing broader context from manifesting even by osmosis."

Thus, mathematics is often mistakenly considered to be value-free among mathematicians as well as parents, employers and teachers. Mathematics teachers often do not believe they are teaching any values at all, rarely consider values in discussions about mathematics teaching, and there is very little research on the subject (Bishop & Clarkson, 1998).

Bishop and Clarkson (1998) suggest three categories of values that mathematics teachers convey through teaching:

- The *general educational* (for example, telling a student off for cheating reveals values of honesty or integrity). "General educational values are qualities which teachers, schools and/or the society/culture aim to inculcate in their pupils, but which are not mathematical in nature. These often have a moral overtone and are essential for the maintenance and enhancement of the social fabric." (Bishop et al., 2001, p. 5)
- The *mathematical* (for example, asking students to describe and compare Pythagorean Theorem proofs reveals values of rationalism and openness). "Mathematical values are associated with the nature of mathematical knowledge itself, and are derived from the way mathematicians of different cultures have developed the discipline of mathematics." (Bishop et al., 2001, p. 5)
- The *mathematics educational* (for example, asking students to show all working in their answers reveals values of mathematical efficiency and effective communication). These are "the norms and practices of doing school mathematics as advocated by mathematics teachers, textbooks and to a lesser degree, perhaps, the school ethos." (Bishop et al., 2001, p. 5)

Bishop and Clarkson present a powerful argument that making these values more explicit and helping mathematics teachers to understand how such values affect teaching generally (regardless of what they may be) improves the quality of mathematics teaching.

## Mathematics, statistics and modelling

Statistics and statistics education have usually been perceived as easier to integrate with ethical thinking than other areas of mathematics because statistical analysis can be thought of as "the study of numbers in context" (Macey & Hornby, 2018, p. 6), with the data themselves being the outcome of some process of measurement. Such measurements always contain a political dimension, reflecting the values of the measurer as much as the technical and practical realities of the data collection process itself (Ridgway, 2016). It is generally accepted that the use of data to model the real world is deeply intertwined with the ethics of what data is collected, how it is collected, and what it is being used to say. Making judgements about what can be inferred, described, or predicted from these data is fundamental to statistical analysis (Velleman, 2008).

This ethical dimension is being brought into ever sharper focus as statistical techniques continue to evolve, increasingly taking advantage of the availability of huge computing power and with ever more readily available and highly-complex data sets. The emerging discipline of data science has the potential to move the learner further away from small, highly curated data sets from which limited conclusions can be drawn, and towards the use of 'open data' (Ridgway, 2016) - large, multivariate data sets collected by third parties such as government institutions, private businesses, and researcher communities. In the world of big data - in classrooms and beyond - consideration of data ethics becomes a fundamental part of the analysis process (François et al., 2020).

Research on good quality learning and teaching of mathematical modelling suggests many rich parallels and overlaps with statistical thinking, for example: cyclical patterns in student progression, encountering and managing mathematical uncertainty, and expressing developing ideas with physical materials and representations (Brady et al., 2020). Below we consider the suggestion that mathematical modelling is not only central to all of school mathematics, but, however construed or constituted, has ethics embedded deep within its structures.

There is good consensus in the literature regarding the relevance and importance of modelling in mathematics education (e.g., Kaiser, 2020). Further, the possibility that modelling is universal, and hence can be applied to any mathematical idea (or system) that students learn about has been suggested. For example, Gravemeijer and Doorman (1999) suggest that context problems, symbolising and modelling are tightly interwoven as students move from 'models of' to 'models for'. "The shift from model/of to model/for concurs with a shift in the way the student thinks about the model, from models that derive their meaning from the modeled context situation, to thinking about mathematical relations" (ibid, p. 119). They describe important efforts to bridge "the gap between the island of formal mathematics and the mainland of real human experience" (ibid, p. 114), by grounding mathematical ideas in situations that feel relevant and authentic to students. They suggest, in particular, the *Realistic Mathematics Education* movement, which rests on the principles of

- **reinvention:** well-chosen, relevant and authentic context problems offer opportunities for the students to develop informal, highly context-specific solution strategies
- **mathematisation**: students' informal solution procedures may function as foothold inventions, or as catalysts for curtailment, formalisation or generalisation.

This suggests one way in which modelling could form part of teaching any mathematical topic. However, some researchers have gone further, stating: "nearly all questions and problems in mathematics education, that is questions and problems concerning human learning and the teaching of mathematics, influence and are influenced by relations between mathematics and some aspects of the real world" (Blum et al., 2007, p. xii).

The question remains as to whether modelling is a useful way of seeing mathematics as part of teaching it, inherently constituted in 'the mathematics' itself, or some combination (or none of these). Brown and Ikeda (2019) suggest many curriculum and resource designers need to catch up to research in this regard, citing Germany as an example of good practice in terms of integration of modelling within national standards across age groups.

Mathematical modelling has been defined in many ways, usually with a focus on both the *use of* mathematics and on the *grounding of* it in some context, perhaps referred to as the 'real world' (e.g., Arseven, 2015), a term that the authors find problematic as it creates an artificial separation between the perceived non-real mathematics and its impact. Brown and Ikeda (2019, p. 233) suggest that "modelling occurs when teachers, students, mathematicians, and others attempt to describe some aspect of the real-world in mathematical terms in order to understand something better or take or recommend actions", a definition which enmeshes ethics at its core with that word 'actions'. There are also multiple researchers that explore ideas of sense-making, sense-checking and sensibility of solutions (e.g. Stillman, Kaiser & Lampen, 2022). In modelling, "any solution that does not make sense in the real-world is no solution at all" (Brown & Ikeda, 2019, p. 234). Modelling also strongly pulls into focus ideas of *authenticity* of task design in mathematics education - questions of whether contexts are incidental, minor and surface level; or embedded, major and deeply related to the problem to be solved (Brown, 2019). We further argue that both modelling in mathematics and mathematical models are essentially meaningless without context and interpretation: *only* a human can bridge the gap between the formal symbols and the contextual and experienced world.

Modelling can also involve descriptive models, or prescriptive models - tools for human decision making - also known as optimisation or normative models (Brown & Ikeda, 2019), which are inherently ethical in nature, prompting important questions such as *best for whom?, best according to what criterion?* and *normal for whom?.* These questions are not only ethical, but often social and political, underpinned as they are by questions of values, responsibility, justice and equity. McKelvey and Neves (2021) further argue that optimisation presents itself in three interacting forms: as a mathematical technique using quantification, de-situation, and formalism, as a tool for legitimation to answer social and other problems, and as a social practice building on aspirations and desires. Its success as a (mathematical) technique fundamentally "relies on longstanding colonial and scientific knowledges that apprehend self and social determination through the lens of development, progress, innovation, perfectionism, and so on" (McKelvey and Neves 2021, p. 97). And indeed, in certain situations, it can already be mathematically proven that optimisation problems in AI tend to have unethical solutions unless the ethics is specifically factored into the objective from the beginning (Beale et al., 2020).

# Ethics in mathematics education

The awareness of the embedded nature of ethics in statistics deviates strongly from the other mathematical sciences. In this section we outline the necessary foundations to incorporate ethics into mathematics education.

**Social responsibility and mathematics education**

Many have argued that the foundation of a just and democratic society is both equitable access to and 'good' use of mathematics. D'Ambrosio (1990, p. 21) suggested, even some thirty years ago:

"Science and technology depend on mathematics, our economic systems are analyzed and regulated through mathematical methods, our opinions rely on statistical data, our health is controlled by the indices of our body chemistry, we are socially classified according to our income, our daily life is paced by schedules, dates, and times. Indeed, our life is regulated by mathematical indices. [Mathematics education] should be looked upon as something that prepares for full citizenship, for the exercise of all the rights and the performance of all the duties associated with citizenship in a critical and conscious way."

Müller (2022) suggests that despite some differences, the two most prominent ways to deal with ethical problems surrounding mathematics, typically known as 'Mathematics for Social Justice' and 'Ethics in Mathematics', are complementary rather than competing educational strategies due to their shared activity of critical questioning and raising social awareness, shared aims of empowering students, and shared epistemological underpinnings of challenging what it is to do and know mathematics. These shared principles form part of the foundation of our proposed framework.

# Ethnomathematics

D'Ambrosio (1990, p. 23) defines ethnomathematics as the idea that mathematical knowledge is inherently a socio-cultural type of knowledge: accumulated through personal experience with the environment and experienced through power dynamics. "Every single manifestation of art, production, behavior, leisure, planning, and so on all the activities of a human being - brings with it numbers, figures, symmetries, harmonies, regularities, extent, design, logic, and so on, and all of this builds up into ethnomathematics." The implication, if we accept this view of mathematics, is a reimagining of teaching mathematics not just as transmitting information, but also creating a just and democratic society by modelling just and democratic behaviours at every level.

Many researchers acknowledge that this tension between the cold, hard stereotype of mathematics they have received and the rather more messy subject as conceived here is difficult for some to swallow. Pais (2011, p. 210) even goes so far as to use the term 'sully':

"'Ethno' shifts mathematics from the places where it has been erected and glorified (universities and schools) and spreads it to the world of people, in their diverse cultures and everyday activities. An ethnomathematical program sullies mathematics with the human factor: not an abstract human, but a human situated in space and time that implies different knowledge and different practices."

Here, we see a shift (conceptualised as a physical move as much a metaphorical one) from mathematics as elite, pedestalised and gatekept to familiar, democratised, and accessible. Both the criticism and the revolutionary cry is that ethnomathematics is letting mathematics out of the ivory tower; letting it, in fact, run riot. This is only sensible, since if "it is now evident that one can wield practically *all* branches of mathematics both for good and harm" as Chiodo and Müller (2018) argue, and if it also affects all of us, then unnecessary participatory hurdles and excessive gatekeeping contradict the visions of democracy based on both shared values and shared access to mathematical thinking proposed earlier. In the context of ethics in mathematics, we, therefore, must go beyond El-Mafaalani's succinctly put words "educational disadvantage is life disadvantage" (El-Mafaalani 2020, p. 95), and say that overcoming educational disadvantage is an absolute necessity to participate in our mathematised world. Consequently, we must not only put lessons from ethnomathematics, but more general issues of equity and inclusion in mathematics learning at the forefront of solving the problem of seeing and responding to the ethics that it is in mathematics.

# Equity

Guterriez (2002; p. 153, in 2012) defines equity in mathematics education as "the inability to predict mathematics achievement and participation based solely on student characteristics such as race, class, ethnicity, sex, beliefs, and proficiency in the dominant language." At first glance, there is an apparent argument that context-free, 'pure' mathematics without embedded ethics may appear more equitable. It certainly feels cleaner, easier, and we may convince ourselves that all students may have a meritocratically equal shot at attaining highly in these sorts of lessons. However, it has convincingly been argued that this is simply measuring "how well students can play the game called mathematics" (Gutierrez, 2012, p. 20), 'achieving' against standards that privilege White, male, Western, 'rational', lone, written ways of working with mathematics that we have mistakenly encoded as neutral (Stinson, 2013). As Ernest (2021, p. 3144) suggests, "there is no acknowledgement of the irony that adopting the position that mathematics is largely value free is indeed itself a values position, or at the very least a meta-values position, leading to the 'choosing not to choose' fallacy". Indeed, some of the authors have directly encountered this fallacy in action. For instance, mathematicians who took a strong position on mathematical modelling - including the modelling done at the beginning of the COVID-19 pandemic - and who insisted that there was "no ethics" in their decision making, with all their choices and outcomes guided by the application of mathematical models, and hence simply "the correct thing to do" (Chiodo & Müller, 2020).

Ethics, then, is embedded in mathematics education through values, and those values can be explicit or implicit, thoughtful or thoughtless, invisible or reified - and they can support or deny equity efforts, both through intent and through neglect, acts of commission or omission. "Equitable classrooms are reflections of a [larger] pedagogical, political, and moral vision" (Lotan, 2006, p. 526), and paying attention to the ways that equity may be expressed in the mathematics classroom is the first step towards this vision for a mathematics teacher. To achieve this, Gutierrez proposes consideration of four dimensions of equity across two axes: Access and Achievement, and Identity and Power, where the notion of Power, in particular, relates to ethics in mathematics education in these ways:

- Who has a voice in the classroom (and who is silenced)
- Students using mathematics to critique society
- Rethinking mathematics, and knowledges more generally, as human endeavours and socially embedded

Embedding explicit ethics teaching in mathematics education is therefore compatible with Gutierrez' conception of equity, and also sits aligned alongside their findings of what it takes in terms of pedagogy to support marginalised students in mathematics: knowing students in meaningful ways that build upon their cultural and linguistic resources; scaffolding learning onto learners' previous experiences without 'watering down' the curriculum; creating classroom environments that are collaborative, friendly, and make good use of group work.

## Ethics in statistics education

The development of the discipline of statistics itself is tangled up with questions of ethics, with many of the much-celebrated historical pioneers of statistical techniques developing the field in the service of racist ideas and ideologies (Rycroft-Smith & Macey, 2022a). More recently, some statisticians' professional associations have sought to embed ethical values throughout the discipline. The American Statistical Association, for example, identifies "promoting sound statistical practice to inform public policy and improve human welfare" as one of their goals (About ASA, 2022). Therefore, there have been growing calls within statistics education to consider issues of both ethics and social justice as important elements of statistical literacy (e.g., Lesser, 2007; François et al., 2020; Bargagliotti et al., 2020).

Much recent research into statistics education advocates for approaches to teaching and learning that involve problem-solving and informal inference within the context of data modelling - posing meaningful questions and using data and representations of data to resolve them (Ben-Zvi & Makar, 2018). Thus. embedding ethics into statistics education means embedding ethics at every stage of statistical activity. Advocates of a problem-solving approach recommend applying a data modelling cycle which may take numerous forms; one example is PPDAC (Ben-Zvi & Makar, 2018), which refers to Problem, Plan, Data, Analysis, Conclusions. The *Pre-K-12 Guidelines for Assessment and Instruction in Statistics Education II (GAISE II)* (Bargagliotti et al., 2020) use a four-stage model comprising: formulate statistical investigative questions, collect data/consider data, analyse the data, and interpret results. Table 1 gives some examples of ethical considerations at each stage:

| **Formulate statistical investigative questions** | Why is this question being asked? How might the answers to the question be used? Who is asking the question, and what is their relationship to the data sources? |
|---|---|

|  | What assumptions are the question based on? |
|---|---|
| **Collect data/consider data** | What does it mean to own data?<br>Who does the data belong to?<br>Which sources of bias influence the data?<br>What is being measured and how has this measure been derived?<br>What is not being measured?<br>How might measuring this quantity influence what is being measured? |
| **Analyse the data** | Why are particular statistical techniques being used?<br>Is a particular technique appropriate?<br>Which features of the data are highlighted or hidden by different choices of representation?<br>How are patterns enhanced or diminished by choices of representation? |
| **Interpret results** | What biases might be reflected in the conclusions?<br>What assumptions have been made in drawing any conclusions?<br>Are there alternative conclusions that could be drawn?<br>Who might be harmed or empowered by the conclusions? |

**Table 1: Examples of ethical considerations at each stage of the GAISE II framework**

A further recommendation suggests the use of large data sets and open data (e.g. Ridgway 2016; François et al., 2020). The use of such data sources create a need to consider contemporary issues of data use, including but not limited to: non-discrimination, fairness, privacy, data protection, and data security. Ridgeway, for example, points out that anonymisation of data is complex and challenging, when seemingly anonymised data can be sliced in such a way that, when combined with additional information, participants can be readily identified. Rycroft-Smith & Macey (2022b, p.3) also suggest some design guidelines for mathematics education assessment and resources which explicitly consider use of data in the classroom, suggesting, for example, that designers should "explicitly encourage students to critique the binary categorisation of sex as a unit of analysis" and "avoid questions with harmful superficial associations, in particular where the meaningful statistical work needed to investigate the deeper issues may be beyond students' current capabilities, or insufficient time is given to exploring the issues in detail; for example, 'associations' between IQ and 'race'." These issues and others are rich seams to explore in the ethical statistics classroom.

It is well documented that statistics can be used for harm, and particularly relevant to modelling is the process of *mathwashing*, which has variously been described as deceiving people about the objectivity of mathematical models with baked-in biases (Mok, 2017), as "a way of lending credibility to an argument or justifying a course of action by appealing to mathematical authority" (Rycroft-Smith & Macey, 2020), and as hiding "a subjective reality under a thin layer of 'objective' mathematical processing" (Allo, 2018, p. 2). The term was originally coined by Benenson in an interview (Woods, 2016) "to describe [...] using the objective connotations of mathematical language to describe products and features that are probably more subjective than their users might think", a particularly useful concept when discussing mathematical ethics. Further, Aragão and Linsi (2022) outlined four strategies through which statistics can be intentionally misused: by outright manipulation, (politically) motivated guesstimating, opportunistic use of methodology space and indicators-management (manipulation) through indirect means. It is clear that teaching ethics as part of statistics should explore not only how to make good decisions and to do good, but also how to avoid harmful decisions and do harm, as well as recognise when harm is being perpetuated and use one's voice to speak up about it. This is the argument we extend to all of mathematics education, as we elaborate further.

# Why is such a framework needed?

Chiodo & Bursill-Hall (2019, pp. 5-6) argue that those in the mathematical community are encouraged to see mathematics as beyond ethics (Platonism), and that, in contrast to law and medicine which teach profession-specific ethics, social and ethical consequences are seen as just not a mathematician's problem (exceptionalism), prompting an urgent need for integrated teaching of ethics in mathematics. The described position to see mathematics as beyond ethics is closely aligned with Kant's position of mathematics as synthetic a priori, universal knowledge gained through a logically reasoned synthesis of axioms and previously proven results. These purely rationalist - and externalist - definitions of mathematics ignore mathematicians and their values, and can therefore often hinder mathematicians and mathematics educators to perceive mathematics in a wider context (Müller, 2018). Definitions like these may also prevent them from understanding ethics in mathematics in the first place, since such ethics always attends to mathematics in context. But as Gutierrez (2012) reminds us, all good teachers focus on context and recognise the positioning of students, problem, mathematics and classroom as crucial to the availability of resources that the student may draw on to express their mathematical thinking and be understood.

It is worth pointing out that postponing the ethical education of students until university is not an option because "students need to train their ethical reasoning just like they train mathematical reasoning" (Chiodo & Bursill-Hall, 2019, p. 40). Ethical reasoning in students needs to be fostered and develops more deeply over time. The ethics of mathematics education is multidimensional and its ethical practice is necessarily dialogical in its nature (Boylan, 2016), but when students learn at school that there is no ethics in mathematics or do not learn about its existence, then it generally has two direct consequences:

1. At university, the lecturers need to try and "undo" that teaching which is (speaking from experience) very hard to do, and
2. University students tend to (from our experience) engage less, or not at all, with ethics in mathematics training, courses, and events, because they just don't see how it could exist or be an important thing to pay attention to.

Teaching ethics in mathematics at universities - both for students majoring in mathematics and those studying mathematics as part of another scientific or professional degree - is already a difficult undertaking without the added burden of having to undo previous harmful lessons which trivialised or even dismissed ethics in mathematics. We found that ethical perseverance, i.e. a "willingness and consciousness of the need to pursue ethical insights and truths despite difficulties, obstacles, and frustrations" (Paul & Elder, 2006, p. 35), and ethical courage, i.e. "the willingness to face and assess fairly ethical ideas, beliefs or viewpoints to which we have not given serious hearing, regardless of our strong negative reaction to them" (Paul & Elder, 2006, p. 35), are particularly hard to teach if students had no previous exposure to ethics in mathematics.

We feel that this makes a strong case for ethics in mathematics being taught in schools, across all mathematical topics, and in all age groups.

## Teacher education in mathematics

Building on the works of Stein & Kaufman (2010), Xenofontos (2016) and Gutiérrez (2013), Xenofontos argues that "[t]eachers stand between intended curricula/policies and actual learning outcomes, mediating mathematics learning through particular knowledge, beliefs, ideologies, identities, experiences, and instructional practices." (Xenofontos et al., 2021, p. 136). Teachers can be supported in examining their views of mathematics as certain, Platonic or infallible by persuasive arguments about productive ambiguity leading to creativity in both teaching and learning (Marmur & Zazkis, 2022). Similarly, examples of professional development programmes with a view to increasing mathematics teachers' understanding of and engagement with social justice include:

- Meaney et al. (2009) who examined teachers' sense of belonging and identity in a Maori-medium school, drawing on Wenger's modes of belonging
- Ndlovu (2011) who developed school-based professional development using theories of reflexive competence in school-university partnerships in South Africa, considering 'socially just' models of improvement for mathematics and science teachers
- Wright (2016) who developed professional development for teachers in England based on principles of: collaborative, discursive, problem-solving, and problem-posing pedagogies; recognizing and drawing upon learners' real-life experiences in order to emphasize the cultural relevance of mathematics; promoting mathematical inquiries that enable learners to develop greater understanding of their social, cultural, political, and economic situations; facilitating mathematical investigations that develop learners' agency, enabling them to take part in social action and realize their foregrounds; and developing a critical understanding of the nature of mathematics and its position and status within education and society.

All these programmes have in common the aim of supporting teachers in understanding their identity and therefore position in relation to mathematics and society. Further, although the terms morality and ethics are often used interchangeably, they are philosophically distinct; "morality is [a] first-order set of beliefs and practices about how to live a good life, ethics is a second-order, conscious reflection on the adequacy of our moral beliefs" (Gülcan, 2015, p. 2624). Principles of effective teacher professional development design suggest it is beneficial for teachers to reflect on their own ideas (Körkkö et al., 2016), beliefs (Maas et al., 2017), and even identities (Beauchamp & Thomas, 2009). We believe there are therefore likely to be benefits to encouraging teachers to reflect on:

- their positioning as mathematicians within society
- their dynamic identities in relation to mathematics
- their experiences of using mathematics in the classroom and in education more widely with ethical consequences

There are many interesting models of teacher professional learning that are consistent with, and would be supported by, our framework; for example the zoom in and zoom out approach to mathematics teacher education (ZIZO) which aims to stimulate the development of a problematising stance toward knowledge production and teaching of mathematics through exploration of historical sources and collective discussion (Moustapha-Corrêa et al, 2022).

As we have previously suggested, in this paper we will take a pragmatic view on ethics. Our framework will be sufficiently general that individual teachers and educators can use it independently of their own ethical theories, meaning that no matter what ethical theories the readers ascribe to themselves, or previous (in)experience they may have in working with ethical dilemmas, our proposed framework will be useful for them.[5] [6]

## Making the case

Danks (2021, p. 3) carefully scrutinises the so-called *neutrality thesis* for digital technologies and finds that not only do digital systems sometimes make ethical decisions themselves, but "the creation of digital technologies involves an enormous number of ethical choices that thereby embed or implement values in the technology itself. Most obviously, digital technologies are almost always designed or optimized for success at a specific task, but 'success' is an ethically substantive term."

---

[5] We found Michael Sandel's *Justice* course to be an accessible introduction to specific ethical theories. Segments from the lectures and the accompanying discussion guides can also be used in a classroom setting to explore selected ethical theories on a deeper level. They are easily adapted to a more mathematical audience. http://justiceharvard.org/justicecourse/

[6] We refer to Blackburn (2021) for a general introduction to ethics, Martens (2022) for data science ethics, Held (2006) for the ethics of care, Briggle and Mitcham (2012) for an introduction to the ethics of the (natural) sciences, and to Crawford (2021) to learn more about the planetary costs of artificial intelligence.

In mathematics education, digital technologies are not only employed *in service of* teaching and learning, but there is considerable overlap between the technologies themselves (or algorithms, or systems, or models) and what may be understood to constitute the mathematics curriculum, such that the second argument applied by Danks is clearly equally applicable to school mathematics - both in terms of content but also, interestingly, in terms of measuring success. For example, when *students* learn about calculating averages, they learn to use simple algorithms that nevertheless contain embedded privileging of certain ideas (norms, accounting for outliers, regression, reducing deviance). At the next level up, when *teachers* use averages to report on student data, they are also communicating embedded privileging of further ideas based on those concepts – just who does the 'norm' and the 'deviant' represent in terms of the assessment criteria and their concomitant assumptions, for example? Likewise, even outside the classroom, the question of *whether* to use an average can have far-reaching consequences, such as the "Robodebt" example in Australia where the government chose to use averaging from annual income to determine who had been overpaid with welfare benefits, even though eligibility was legally determined by fortnightly income and not annual income (Henriques-Gomes, 2020). In this instance people were sent (unjust) debt recovery notices in the tens of thousands of dollars, and there were reports of people becoming the victims of suicide as a result (McPherson, 2019). As Danks suggested earlier, success is indeed ethically substantive in definition, and as mentioned earlier concepts of optimisation which underlie models always include embedded values, however deeply hidden.

As the text-based mathematical problems that students encounter during their school education have traditionally been set in a context familiar to them, we are bound to see more textbooks containing problems involving recent technologies and social issues. Many of these text-based problems already include elements of optimisation, and regularly fail to mention the wider (ethical) context of what is being optimised. Shulman (2002) exposed some particularly dangerous examples from textbooks, including a question on poisoning a city using gas and the problem of setting up an offshore oil well without taking environmental factors into account. Questions of success/optimality, and hence ethics, are therefore not only present in the student teacher relationship and the education system as a whole, but in the mathematical problems that students are being asked to solve.

The general argument for teaching some kind of ethics embedded in mathematics at school level is also not new, although its urgency is added to by rapid technological development. In a conference paper Bishop et al. (2001, p. 9) ask:

"Could it be accepted that one would discuss "the social responsibility of mathematicians"? Why does that seem such a strange idea? After all mathematics permeates all our lives and will increasingly do so. So why should a good, modern, mathematical education not face students and teachers with the moral, ethical, and social issues of applying mathematics everywhere?. Teachers, learners, and everyone who desires an education fitting students for their roles in a more democratic society needs to understand about the role that values play in a mathematical education."

# Useful case studies

While there are many case studies to be found in the literature (e.g., O'Neil 2016), below we give two example areas: financial mathematics, and bus timetabling as an application of AI and optimisation. The most common ethical issues often stem from basing decisions on poorly founded or badly articulated assumptions that later cause problems - issues to which we can all relate. Often in mathematics, we produce semi-predictive tools that we hope tell us something useful about the future behaviour of an object or process. But quite regularly these tools are not as predictive as we wish and instead come with severe limitations. As mathematicians, we know this, but sometimes either:

1. Fail to communicate the limitations effectively to others, or
2. Fail to properly acknowledge and stay within the limitations ourselves.

We now explore these with some case studies on mathematical work and its impact.

**Financial mathematics**

One area where we have seen problems like this arise is mathematical finance, where mathematics is used to give predictions regarding future price movements. The failure to adequately acknowledge mathematical limitations and mathematics' "precarious abstraction" has played a critical role in many large-scale issues, including the Global Financial Crisis of 2008-2009 (MacKenzie 2011, 2012), when mathematicians "enabled bad debts to be hidden by devising Collateralised Debt Obligations (CDOs) that bundled together thousands of mortgages and sold them off in slices with precisely calculated risks attached" (Mendick, 2020, p. 103). Interesting explanations of these ideas can be found in the films *Margin Call* and *The Big Short*, and Mendick (ibid) provides a useful analysis of the presentation of questions of ethical mathematics in these films in *The financial crisis, popular culture and maths education*.

But in today's financial world, mathematics is not just used to model markets and to price options and other derivatives. Instead, through high-frequency trading, mathematics sits at the core of most trading activities when the algorithms interact at speed beyond human comprehension (MacKenzie 2015, 2016, 2018). To put it bluntly: modern finance only exists because of mathematics, and when it goes wrong, it can go very wrong.

**Designing bus routes**

It is not just mathematical finance where mathematics can cause problems. Many mathematicians go on to work in technical and seemingly 'neutral' areas. The skills we learn as mathematicians - such as ordering, classifying, counting, measuring, sorting, stacking, logical reasoning, and problem-solving - make us perfect candidates to work in and contribute to many other fields. You may already be familiar with problems in AI around "killer robots" and racially biased facial recognition algorithms. But what about something much more innocent. What about bus timetabling? We can, after all, use AI and optimisation techniques to help us design a bus network.

There are many steps to solving the problem, and each of which can raise serious issues if not treated carefully, as we show in Figure 1.

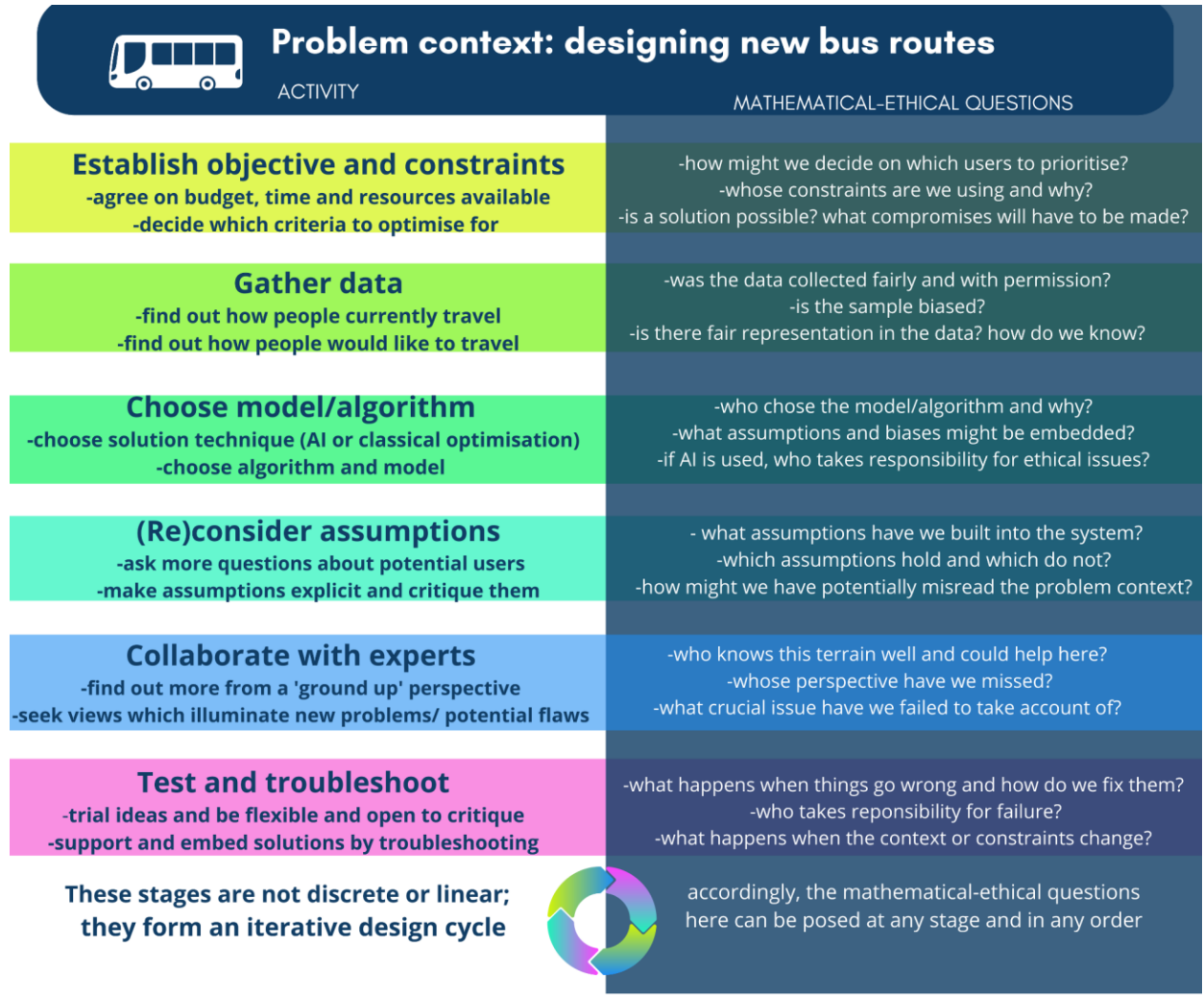

Figure 1: a mathematical-ethical exploration of designing new bus routes

As mathematicians working on this project, we may be tempted to want to 'keep our hands clean' and only perform computations, but the mathematical work is entwined with work that is also ethical, social, and political. It may turn out that children can't get to school, or doctors can't get to work, or traffic gridlocks in part of the city, or shops go out of business due to the change in footfall. In this context, it must be noted that older data is rarely suitable to address these issues as the needs of women and caretakers taking public transport are often not properly reflected in it (Perez 2019).

If there are citizen protests in the streets because of the perceived failure of our new transport system, we can't just tell people to "please wait while we re-run all the computations"; the time for that will have passed. Nor can we tell them that "we are just the mathematicians; it's not our problem" because it very much is our problem; a problem we created, and one that we must fix.

Thus, we often need to go through the optimisation cycle many times before we can start considering deploying our solution.

All of the problems we have outlined in these case studies apply much more widely. They are ubiquitous in AI and machine learning, issues of cryptography and surveillance, statistics, mathematical biology, algorithmic modelling, climate modelling, and many more areas of mathematics. Particularly, the ethical issues in AI and algorithms are now widely discussed in the literature, but we hope that the reader now sees that they can be found throughout mathematics. What matters most here is not the specifics of the particular area of mathematics, but rather the general approach that mathematicians take to their work and to problem-solving. Even something as simple as "What will you do when you make a mistake?" (which can happen in any area) is occasionally replied to by mathematicians with "If I've done it properly, I won't make a mistake (in my mathematics)", and therefore they overlook the importance of (non-mathematical) post-deployment response strategies. If you honestly don't believe you could do something wrong, why have a plan for when something goes wrong? This, and many other issues, transcend all mathematical work. Having non-mathematical response strategies has proven to be particularly important in the modelling of Covid-19 (Chiodo & Müller 2020).

Having seen that ethical questions arise at all steps of mathematical modelling, one could even begin to move discussions in the mathematics classroom to a more abstract level by asking, "If a mathematical process makes an ethical decision that causes problems, who do we hold accountable?" When students realise that you cannot put a computer in jail, they quickly trace back, and it strikes them that the ethical accountability rests with the creators.

# Framework

We draw on the Chiodo and Bursill-Hall (2018) level-based system for mathematicians' ethical engagement or consciousness, suggesting these levels are neither discrete nor linear exactly, but a useful framework from which to begin more nuanced discussions of this important topic that will broaden and deepen over time, and translating the levels into use for mathematics teachers in the following ways:

| Chiodo & Bursill-Hall : ethics for mathematicians levels | Translated for mathematics teachers… | So that mathematics learners... |
|---|---|---|
| Level -1: Actively obstructing efforts to address ethics in mathematics | Actively obstructing efforts to address ethics in their mathematics classroom | Are encouraged to dissuade peers from thinking about the social impact of mathematics, and the values encoded in mathematics and its education |
| Level 0: Believing there is no ethics in mathematics | Believing there is no ethics in mathematics education (and mathematics) | Believe mathematics is pure, neutral, apolitical, value-free, cannot be used for harm (and therefore are more likely to use it for harm, whether through malice or incompetence) |
| Level 1: Realising there are ethical issues inherent in mathematics | Realising there are ethical issues inherent in mathematics education (and mathematics) | May not receive any, or very little, discussion or discourse around ethics in mathematics education |
| Level 2: Doing something: speaking out to other mathematicians | Doing something: designing or adopting questions, tasks, or lessons that address ethical issues in mathematics education | Begin to experience discussion or discourse around ethics in mathematics education – which may or may not be good quality, or embedded in classroom experience that is aligned/coherent |
| Level 3: Taking a seat at the tables of power | Thoughtfully posing and exposing ethical dilemmas in the classroom; developing *political conocimiento* (Gutierrez, 2017); taking ethics in mathematics beyond their own classroom and beginning to also bring it to other contexts | Experience rich, structured discussion or discourse around ethics in mathematics education and are safe and belong in an environment where these ideas are aligned with the classroom climate |
| Level 4: Calling out the bad mathematics of others | Using ethical awareness to object to bad data in education systems; questioning and contributing to policies; exposing mathwashing; auditing and challenging with mathematics; confronting gatekeeping | Have access to mathematics as a tool for empowerment, justice and belonging<br>See everyone by default as a mathematician<br>Are aware of the embedded ethical decision points at every level of problem solving |

**Table 2: Ethics framework for mathematics teachers**

Below, we discuss each level in more detail, exemplifying what it may look like in practice for mathematics teachers.

## Level -1: Actively obstructing efforts to address ethics in mathematics

At this level, mathematics teachers use their **voice** to suppress the idea that mathematics and ethics are intertwined, and to **silence** students and peers that suggest that ethics have an important place in mathematics education. They may make use of their positions of power to intimidate and use arguments based on authority and status. They publicly draw on flawed arguments of the exceptionalism, neutrality or value-free nature of mathematics, using discourse such as:

- "Those sorts of lessons belong in other subjects, not in mathematics."
- "Stop asking silly questions and focus on the real mathematics."
- "Do you actually think algebra can kill people?"
- "Stop misbehaving and get back to the mathematics lesson!"
- "That's not real mathematics."

Many of these discourses may not always happen consciously or stem from bad intentions but instead as the result of received wisdom or a pyramid scheme of math abuse (Fiore, 1999). Speaking negatively about ethics in mathematics can, in our experience, occasionally overlap with other problematic pedagogical approaches that may support deficit narratives around minoritised identities, perpetuate gender stereotypes in the mathematics classroom, fail to examine the teacher's own biases, and/or perpetuate inequities through gatekeeping and withholding access to particular, privileged mathematical 'ways of knowing'.

Sometimes these discussions are also improperly closed off for the mere fact that they can attract controversial and polemical statements, because they are unpopular with a small number of very loud voices (e.g., a few parents) or because discussion spaces in mathematics are restricted to a very narrow definition of mathematics, and people are unsure where the right place to discuss these is. For an example of such discussion, we refer to the many user comments in the mathoverflow thread on collecting "references for literature from mathematicians who provided critiques and proposals concerning ethical aspects of mathematics research"[7].

## Level 0: Believing there is no ethics in mathematics

At this level, mathematics teachers are largely **silent** on the issue of mathematics and ethics. They do not use their **voice** to question, support or discuss ethical-mathematical issues with students or peers. They largely do not think about these issues in any depth at all, and if they have thought about them, they will generally try to remove the human component from mathematics, focus on mathematics as the only form of universal and necessary knowledge, or fall back on the idea that 'all tools/artifacts are neutral'. They may privately or publicly draw on received, surface or flawed arguments of the exceptionalism, neutrality or value-free nature of mathematics, using discourse such as:

- "I got into mathematics for the right and wrong, black and white – not this stuff."
- "Mathematics isn't about being human – it is transcendental, universal, beautiful."
- "Mathematics is just the tool."
- "Mathematics is truth, and truth is beauty."

---

[7] https://mathoverflow.net/questions/346894/references-for-literature-from-mathematicians-who-provided-critiques-and-proposa/346897#346897

In our own attempts to teach and discourse about ethics in mathematics at universities, we found that mathematicians at this level are largely unaware of, or disagree with, the debate on "do artifacts have politics?", begun by Winner's (1980) essay with the same title, and that a general awareness of the various disputes in the philosophy of science from the 20th century helps to move quickly from level 0 to a deeper understanding of ethics in mathematics.[8] Recent case studies on epistemic injustice in mathematics (Rittberg et al., 2020) and the many accessible works on mathematics for social justice (e.g., Buell and Shulman, 2021) have also proven useful as a first point of contact when speaking to mathematicians and teachers who do not believe that mathematics has ethics. In this context, we particularly want to mention Paul Ernest's (2016b) "Dialogue on Ethics of Mathematics" as an introduction to the topic, which we suggest as valuable for teachers to reflect on and critique as part of their journey into awareness of ethical-mathematical issues.

## Level 1: Realising there are ethical issues inherent in mathematics

At this level, mathematics teachers begin to consider ethical issues in mathematics and mathematics education, perhaps through coming into contact with stories in the media, case studies or interesting resources. They may still consider that ethical issues do not apply to all of the mathematics curriculum. They are likely to remain **silent** at this stage, but crucially are beginning to listen to **voices** that allow space for different narratives around ethical-mathematical issues. They will have stopped actively professing that mathematics is value-free and ethically neutral. They may start to engage in and to notice themes in discourse regarding individual examples of bad mathematics, the ethics of individual subfields of mathematics (e.g., finance, machine learning), and about mathematics as a human activity, using and recognising ideas such as:

- "Mathematics has a very real impact on the world."
- "We shouldn't do harmful things with mathematics."

---

[8] In this context, we particularly want to mention Kuhn's (1962) rejection of science as a value-neutral, universal and impersonal knowledge, Polanyi's (1946, 1958) ideas on authority in science and scientific communities, Ravetz's (1971) ideas on science as a craft and the societal obligations coming from scientific solutions, Lakatos's (1978) work on the messy nature of scientific discovery, Latour's (1987) arguments for a socially constructed nature of scientific knowledge, Jaggar's (1989) feminist epistemology on emotion and values as embedded in science, or, more situated in mathematics education, Walkerdine (1990) on rationality and the construction of the masculine 'reason' as control.

- "I read an interesting article yesterday about the involvement of mathematicians in some scandalous event."
- "Ethics is part of all human work and decision-making, so of course ethics is also part of mathematics."
- "Statistics are often misused, and we should be aware of that."
- "Here is an example of mathematical work causing harm to a particular group of underrepresented/minoritised/vulnerable people."

We have found that mathematicians and teachers at this level can occasionally feel rather lonely in this endeavour, in that they may have a hard time finding peers to engage in the discourse. To such mathematicians and teachers, we would like to say that we have observed many more peers at level 0, who only need to be gently encouraged into the conversation by someone, than people at level -1 who actively try to suppress ideas on ethics in mathematics. To find one "friend" at this level, you need only raise it with a handful of colleagues. And having a friend to talk to makes the endeavour much more productive and enjoyable.

Finally, it is worth noting that talking about ethics in mathematics is also a way to protect the mathematical enterprise as a whole, or as beloved mathematician Andrew Wiles said in the Times (Devlin 2013): "one has to be aware now that mathematics can be misused and that we have to protect its good name."

## Level 2: Doing something: speaking out to other mathematicians

At this level, mathematics teachers begin to tentatively use their **voice** to address ethical-mathematical questions in the classroom. They may design or adopt questions, tasks, or lessons that address ethical issues in mathematics education, reflecting on their experience and beginning to consider the far-reaching consequences of ethics in mathematics and mathematics education at many levels, starting to see patterns of who speaks and is **silent**, who counts and is counted. Teachers at this level may also start to focus on the tools required to do better, share their developed ethical mathematics exercises[9], gain new forms of self-awareness and reflection through the discussion with others, and begin to work for changes in a semi-structured way, focusing on either small-scale effects (e.g., setting up a talk, project week, writing mission statements) or beginning to push for larger changes by slowly forming a group of allies. They may pay particular attention to and engage in discourse such as:

- "Mathematics has such a potential to cause harm that we need a Hippocratic Oath for mathematicians."
- "Our department/team should start up a series of talks on ethics in mathematics."

---

[9] Wares et al. have developed exercises to teach ethical mathematics in the context of Covid-19: https://sites.google.com/view-covid-19teachingmodules and the Cambridge University Ethics in Mathematics project provides a short list of exercises for first year university lectures in applied and pure mathematics: https://www.ethics.maths.cam.ac.uk/course/.

- "I've written some classroom mathematics questions that contain some ideas about ethics and social justice- will you take a look?"
- "I've prepared a mathematics project that my students can do which gets them to explore the ethics of this piece of mathematical work."
- "We've set up a group of interested people - mathematics teachers and students - to speak out about a particular injustice that mathematical work has created."[10]

Mathematicians at this level will also begin to speak about their own experiences and identity, addressing their own (false) assumptions and misunderstandings from the past, as well as the student-teacher relationship in the mathematics classroom.

- "I finally get that the way [the instructor] taught her class WAS about social justice… that teaching mathematics about, or through, social justice isn't just about poverty statistics and world population figures... it's also in the thoughts and actions of the teacher toward his/her students and in the thoughts and actions of students toward each other. It's about feeling safe to be who I am and, at the same time, to critically question who I want to become and what (and who) I value. And, most of all, I think it's also about opening up the content of mathematics [...] to this same kind of critical questioning." (fictional course evaluation of a future teacher, in Nolan 2009, pp. 214–215)

Through their interactions with others, teachers at this level also gain a new perspective on caring for pupils. They understand that good mathematics does not exist without ethics, and so part of caring about students is supporting them to be better at mathematics and giving them the appropriate safe space to grow as mathematicians. But they also understand that a safe space does not exist in a discursive vacuum, that silence is not neutral but is a value choice.

## Level 3: Taking a seat at the tables of power

At this level, mathematics teachers thoughtfully pose and expose ethical dilemmas in the classroom. They begin to develop *political conocimiento:* the knowledge that all teaching is political and value-laden, knowledge that critiques, deconstructs and helps to reinterpret systems to better support all students (Gutierrez, 2017). They see that ethical ideas underpin not only the type of access students are permitted to mathematics but also the related implicit epistemological characteristics of the discipline, and relate issues of identity, power, equity, belonging and justice to mathematical-ethical issues. They will start being more active outside their classroom (for example, in organisations, professional associations, professional development contexts) about getting ethics into the mathematics classroom.

They begin to notice who is **silenced** and who is permitted a **voice** in mathematics education, and start to use their own position or their institutional power to advance the cause of ethics in mathematics by participating in local (e.g., school committees, faculty boards, etc) or national

---

[10] The Just Mathematics Collective is an example for such a group. Their website provides a valuable list of resources for interested students and academic staff: https://www.justmathematicscollective.net/.

tables of power (e.g., learned societies, teacher associations). They may participate in discourse such as:

- "We're setting up a curriculum group to work on integrating ethics into our mathematics teaching."
- "I'm joining this teacher group to help promote ethical and responsible work in mathematics."
- "We're organising a workshop on ethics in mathematics, and are inviting other mathematicians and teachers to come and participate in the discussion."
- "We're campaigning to have a new person appointed to this role because the incumbent is actively blocking any progress on teaching ethics in mathematics."

Teachers at this level understand that it is not enough to only teach their students how *to play the game (of mathematics)* but also how *to change the game* (Gutierrez, 2009), and that to *change the game*, they must find their ways to the tables of power themselves. El-Mafaalani (2018, 2020) asks the question: who eats at the table and who still sits on the floor? In the past, only a small group of people sat at the table, but now an increasingly diverse group of people are finding their way there: women, members of the LGBTQIA+ community, immigrants and other minoritised groups. Taking a seat at the tables of power must always come with the commitment to help those still sitting on the floor - or not yet even in the room where it happens.

## Level 4: Calling out the bad mathematics of others

At this level, mathematics teachers thoughtfully design lessons, tasks, and classroom discussions that explore mathematical-ethical dilemmas across all areas of mathematics, supporting colleagues to do the same. They use their **voice** responsibly and carefully: to critique and question lessons that deal with ethics in mathematics in a surface manner or try to 'solve' issues too neatly or without attending deeply to student agency, equity and justice. They use ethical awareness to object to bad data in education systems; question and contribute to policies; expose mathwashing; audit and challenge with mathematics; and confront gatekeeping. They seek out, identify and try to prevent instances of "bad ethics" by other teachers, both in the way others teach mathematics and the way others suppress ethical thinking in their students.

They see, question, and speak up for those who are **silenced**. They actively contribute to discourse such as:

- "I've approached another mathematics department/school because I see their curriculum contains nothing on ethics in mathematics, and the students are missing out on this important learning."
- "Today we're going to learn about the pitfalls of the various modes of mathematics assessment."
- "I've done a study of the statistics to show that our entrance and teaching processes in mathematics are systematically excluding or suppressing certain groups."

- "I have prepared a presentation for the headteacher on the ways that we could use data in the school more ethically."
- "I've written to publisher X about problematic exercises in one of the mathematics textbooks."
- "We have voted as a department not to use Y publisher any more because we consider that they do not practice business ethically in the education sector."
- "We have designed a set of community evenings inviting parents in to address some of the deficit messaging we have heard perpetuated around certain groups of students being 'naturally bad at mathematics'."
- "We are moving to mixed-attainment classes because we do not want to restrict the expectations of our students"

At this stage, teachers are likely to develop an excellent eye for problematic mathematical problems in textbooks, school curricula and outside of the classroom altogether. More detail about teaching at this level can be found in Shulman's (2002) work *The Teaching of Ethics in Mathematics Classes,* and Chiodo and Bursill-Hall's (2019) article *Teaching Ethics in Mathematics*.

# Application

What would it mean to be 'an ethical mathematics teacher'? We recognise that ethics are not video games (although plenty of video games explore ethics beautifully); the goal here is not necessarily that everyone 'achieve' the highest level possible in the shortest time possible, but something rather more complex. As with all frameworks, perhaps the first and most useful benefit is as a map: for users to examine, understand, compare and finally locate themselves and their pedagogical stance, perhaps also beginning to position the views of others they encounter. Since we draw on teacher reflection models of professional learning, we encourage the kind of critical reflective practice that empowers teachers, respecting teacher agency and autonomy to some important degree. We also believe in 'show, don't tell' models: we suggest supporting teachers by giving them the framework and examples of case studies or ethical-mathematical situations, and stepping back. In particular, we hope this framework can be used in support of exploring the *hidden ethics fallacy* we explored earlier - the mistaken idea that those dealing in mathematics can 'opt out' of taking a position at all, by declaring they have made a choice not to engage with ethics (or politics, or social justice, or equity).

Our exemplification of the levels has explored that mathematics teachers who attend to ethics in the ways we propose do not do so as an add-on or a tick-box exercise - not a one-off worksheet, or reading a single article, or attending an hour's training, although all these things may prove useful contributions to their emergent understanding. Instead, they begin to ask deep questions - about themselves, about education, but crucially, about mathematics. It is not for us to provide these answers, but to be delighted that they are being asked, as they should. Chiodo and Bursill-Hall (2019) point out three key aspects of (effectively) teaching ethics in mathematics:

- The need for **standalone lessons** directly addressing ethics in mathematics.

- The need for ethics to be gently woven in **through all aspects** of mathematics education: in mathematical exercises, in projects, and in assessments.
- The need for colleagues to **refrain from dismissing** the importance of learning about ethics in mathematics, be it explicitly or implicitly.

Chiodo and Bursill-Hall (2019) also suggest that such a teaching plan might start off trying to convince students that *there exist ethical issues in some parts of mathematics*, but by the end help students reach the conclusion that *for all mathematics, it is important to consider the ethical issues contained therein*. Their work is directed at teaching ethics at the undergraduate level, but many of their suggestions and principles also apply when teaching it at school.

# Implications

Our aim in creating this framework is to reveal the myriad of ethical microdecisions both within the domains of mathematics curriculum decision-making and mathematics pedagogy decision-making, and to expose the hidden ethics fallacy as used in mathematics and mathematics education contexts. We see this as of clear benefit to a) teachers, b) students, and c) society.

## Implications for teachers

As Ernest (2021) suggests, the first step towards reducing the harm done by mathematics in educational contexts is acknowledging that such harm exists - what he refers to as 'the dark side of mathematics', in three ways in particular:
- students frequently experience harm in mathematics classrooms
- mathematics education certification serves as a 'critical filter' to prevent equitable access to opportunity
- styles of thinking perpetuated in the Platonic mathematics classroom can lead to damaging ideas of contextwashing, mathwashing, or faux 'neutrality'.

Teachers who have self-identified at Level 0 and in particular Level -1 may have some considerable discomfort. Some first small steps for the teachers who have located their thinking in this space may be:

- to change behaviour and pedagogy in small but important ways, for instance, by being aware that some problems in books may be unethical and to leave them out.
- to consider ethics in mathematics as a standalone or extracurricular activity at first: a project week or even only a project hour is always an option.
- to understand that one does not have to jump straight to level 3 or 4. It is not necessary to do things perfectly. Taking small steps in the right direction is plenty at the beginning
- to be aware that taking *the step from speaking out against ethics in mathematics or from not realising its existence to understanding that there is ethics in mathematics* is a big conceptual leap.

- to understand that there are allies. If teaching the first lesson on ethics in mathematics seems daunting, it might be beneficial to ask colleagues from other subjects (e.g., philosophy, religion) for help.

# Implications for professional development designers, teacher educators, resource and assessment designers

## General mathematical-ethical questions setting

Professional development materials for mathematics teachers often recommend drawing attention to a fixed set of dimensions when considering a particular topic to teach, such as the *Preparing to Teach a Topic framework* in the Open University mathematics education courses which draw on ideas of *accessing core awarenesses* (Mason, 2009).

In a similar way, when designing resources, it may be useful to begin with a set of general, embedded questions that apply to all mathematical areas before considering those that may be more specific, with the goal of accessing core awareness of mathematical-ethical questions as universally applicable. We suggest the below:

- Which mathematical notation and symbol system/s are commonly used and where (according to received wisdom/digging deeper) do they originate?[11]
- Which names or attributions (for problems, formulae or ideas) are used in this topic, and where did they come from?
- In what contexts might this mathematics be used, and who might be impacted?
- What type of word problems are frequently used and why?
- How might this mathematics be used for harm?
- How might this mathematics be used for good?
- Have you heard of any news stories involving this mathematics?
- How is this topic assessed and who might benefit more, and less, from this type of assessment?
- What would a *wrong but useful* (Beveridge and Gale, 2013-2020) answer suggest?

An example of this in use with several possible topics is shown in Figure 2.

---

[11] To get a sense for the richness of the history of mathematical notations, we refer the reader to Cajori (2019).

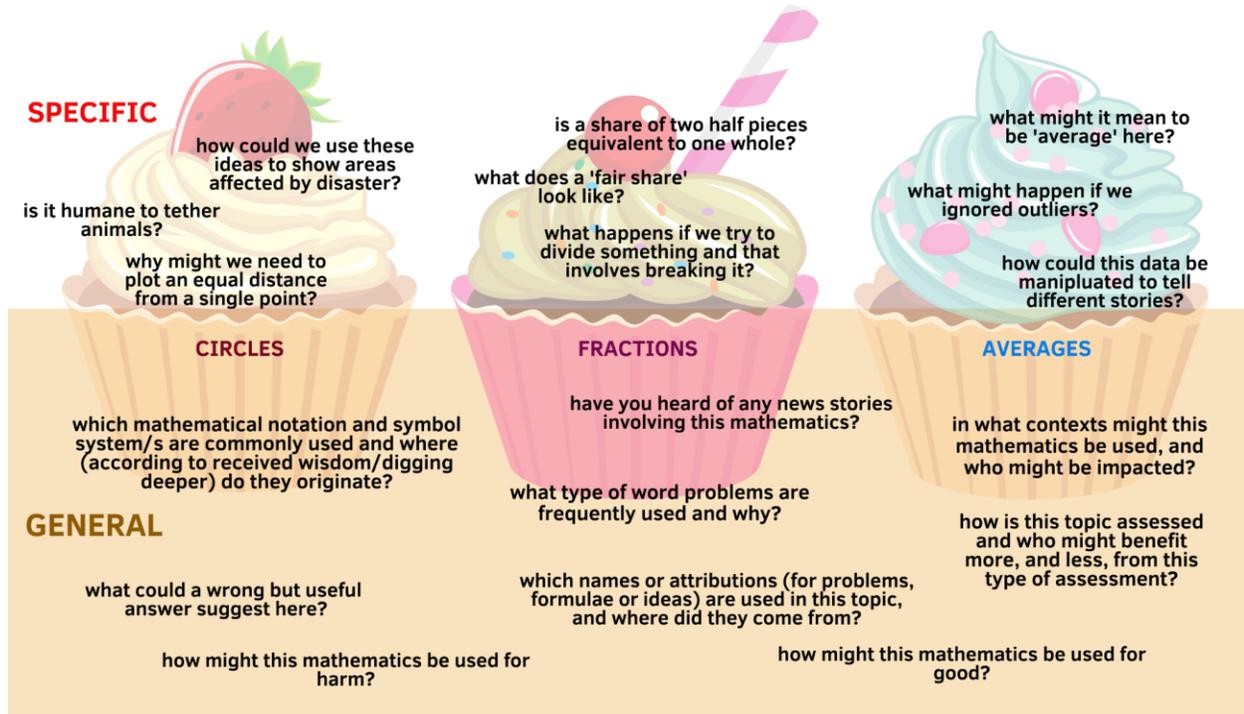

**Figure 2: suggested two-part model (general and specific) for mathematical-ethical design questions**

We also recommend Rycroft-Smith and Macey's (2022b) *Draft guidelines for assessment and resource design in mathematics education* and their discussion of the design principles behind the guidelines (Macey & Rycroft-Smith, in press) and welcome feedback and collaboration.

## Implications for mathematics education researchers

What we attend to and our choice of focus, in research as elsewhere, constructs discourses and the telling of stories - not least with ourselves as an audience. As Gutierrez and Dixon-Roman (2011) suggest, the gaze both confers authority and normalizes. "The kinds of research questions we ask influence the knowledge that is created as well as what we might be able to do with that knowledge" (ibid, p.25). We call on mathematics education researchers to attend to mathematical-ethical questions more widely, to collaborate across disciplines in order to further understand mathematical-ethical issues, and in particular to use and test this framework to destruction, improving on it as necessary.

## Implications for students

It is our wish that this framework is another piece of the puzzle that moves from mere 'gap-gazing' (Gutierrez, 2008) to doing something useful to reframe the problem. Macey and Rycroft-Smith (in press) speak of an *action gap* in education, where it appears all too easy to speak of equity in mathematics in terms of platitudes and saccharine statements, but all too hard to make a start, to actually do something, however small - even though one without the other is nonsense. As teachers and educators, we know that you never, ever 'just' teach mathematics.

"With every single thing about math that I learned came something else. Sometimes I learned more of other things instead of math. I learned to think of fairness, injustices and so forth everywhere I see numbers distorted in the world. Now my mind is opened to so many new things. I'm more independent and aware. I have learned to be strong in every way you can think of it. (Lupe, Grade 8)" (Gutstein, 2003, p. 37)

We hope that the permission and courage to address ethical-mathematical questions - to use ethics and mathematics together to *solve problems* (in the sense of to *make things better*), trickles down to students, and fast; in fact, we hope the floodgates open. We hope students feel empowered to ask for, as is their right, ethics as well as ethical treatment in mathematics lessons.

# Conclusion

Mathematics is an active, dynamic, thriving endeavour, and its purposes and uses are thrillingly diverse. We, as humble mathematicians of varying pedigree, do not seek to render it more static, cold, or abstract; or to freeze it in the glorious past. We know that ethics is the lifeblood of mathematics as it lives and breathes now, and we offer this knowledge to mathematics teachers in the hope that:

- they may feel empowered to explore ethics in mathematics lessons with all students of all ages
- they may feel motivated to explore ideas of *political conocimiento* and what it means in their context
- everyone involved with mathematics, mathematics communication and mathematics education may collaborate in the pursuit of mathematics for justice and belonging and the avoidance of harm.

In this chapter, we have explored the use of mathwashing - using mathematics as a borrowed cloak of legitimacy and credibility by drawing on associations with objectivity. Here, we propose a perpendicular term: *math-smashing*, which we define as using mathematics (and/or mathematics education) to **break** enforced silences, barriers of 'objectivity', and instances of gatekeeping. The metaphor is intentionally destructive, since we see this process as something like breaking a glass ceiling: it is the enabling of access to *what mathematics truly is*, and importantly, we believe that access is *for everyone*. We agree with Mendick (2021, p.106) that "mathematical knowledge is part of the democratic competence we need to hold the powerful to account", and we all need to wield our mathematical knowledge-is-power with intention, care and community.

Mason, J. (2009). Teaching as disciplined enquiry. *Teachers and Teaching: theory and practice*, *15*(2), pp. 205-223, DOI: 10.1080/13540600902875308.

Martens, D. (2022). *Data science ethics. Concepts, techniques and cautionary tales.* Oxford University Press.

McKelvey, F. & Neves, J. (2021): Introduction: optimization and its discontents. In Review of Communication *21*(2), pp. 95–112. DOI: 10.1080/15358593.2021.1936143.

McPherson, Emily (2019, July 30). Queensland man took his own life after learning of Centrelink debt, mum says. *9Nnews.* Available online at https://www.9news.com.au/national/centrelink-robodebts-queensland-man-took-his-own-life-over-debt-mum-says-australia-news/e31e6f28-2e4b-4d3f-9095-d8f74e00cbc1, checked on 20/12/2022.

Meaney, T., Trinick, T., & Fairhall, U. (2009). 'The conference was awesome': Social justice and a mathematics teacher conference. *Journal of Mathematics Teacher Education*, *12*(6), pp. 445–462. https://doi.org/10.1007/s10857-009-9122-3

Mendick, H. (2020). The financial crisis, popular culture and maths education. Chapter 9 in Ineson, G., & Povey, H. (Eds.). *Debates in Mathematics Education (2nd ed.).* Routledge. https://doi.org/10.4324/978042902101

Mok, K. (2017). Mathwashing: How Algorithms Can Hide Gender and Racial Biases. *The New Stack.* Available online at https://thenewstack.io/hidden-gender-racial-biases-algorithms-can-big-deal/, checked on 19/12/2022.

Müller, D. (2018). Is there Ethics in Pure Mathematics? *EiM Discussion Papers* 2/2018.

Müller, D. (2022). *Situating "Ethics in Mathematics" as a Philosophy of Mathematics Ethics Education. Arxiv Preprint,* pp. 1-15. https://arxiv.org/abs/2202.00705

Müller, D., Chiodo, M., & Franklin, J. (2022). A Hippocratic Oath for mathematicians? Mapping the landscape of ethics in mathematics. *Science and Engineering Ethics 28*, 41: Article 41, pp. 1-30. https://doi.org/10.1007/s11948-022-00389-y

Ndlovu, M. C. (2011). University-school partnerships for social justice in mathematics and science education: The case of the SMILES project at IMSTUS. *South African Journal of Education*, *31*(3), pp. 419–433.

Nolan, K. (2009): Mathematics in and through social justice: another misunderstood marriage? *Journal of Mathematics Teacher Education 12*(3), pp. 205–216. DOI: 10.1007/s10857-009-9111-6.

O'Neil, C. (2016): *Weapons of math destruction. How big data increases inequality and threatens democracy.* London: Penguin.

Pais, A. (2011). Criticisms and contradictions of ethnomathematics. *Educational Studies in Mathematics*, *76*(2), pp. 209–230. https://doi.org/10.1007/s10649-010-9289-7
34